# Fixed-point iterative linear inverse solver with extended precision

Zheyuan Zhu, Andrew B. Klein, Guifang Li, and Shuo Pang

*Abstract*— Solving linear systems is a ubiquitous task in science and engineering. Because directly inverting a large-scale linear system can be computationally expensive, iterative algorithms are often used to numerically find the inverse. To accommodate the dynamic range and precision requirements, these iterative algorithms are often carried out on floating-point processing units. Low-precision, fixed-point processors require only a fraction of the energy per operation consumed by their floating-point counterparts, yet their current usages exclude iterative solvers due to the computational errors arising from fixed-point arithmetic. In this work, we show that for a simple iterative algorithm, such as Richardson iteration, using a fixed-point processor can provide the same rate of convergence and achieve high-precision solutions beyond its native precision limit when combined with residual iteration. These results indicate that power-efficient computing platform consisting of analog computing devices can be used to solve a broad range of problems without compromising the speed or precision.

*Index Terms*—fixed-point, iterative solver, analog computing

## I. INTRODUCTION

The problem of solving a large-scale linear system occurs in a broad range of science and engineering problems, such as image reconstruction [1]–[3], compressive sensing [4]–[6], and predictive control systems [7]–[9]. These problems aim to deduce the inputs, **x**, that generate the observations, **y**, given the physical model of a linear system, expressed by matrix **A**. Due to the scale of the system model and noises in the measurement process, direct inverse of the observations **y** is usually impractical and, in many cases, impossible. As a result, many inverse algorithms rely on iterative refinements to approximate the solution [10]. Assuming that the measurements are dominated by Gaussian noise, finding the solution $\hat{\mathbf{x}}$ given the observations **y** amounts to minimizing the objective function [11],

$$\hat{\mathbf{x}} = \underset{\mathbf{x}}{\operatorname{argmin}} \frac{1}{2} |\mathbf{y} - \mathbf{A}\mathbf{x}|_2^2, \qquad (1)$$

where $|\cdot|_2$ denotes the $\ell^2$ norm. The solver starts with an initial guess, $\mathbf{x}_0$, and moves along the gradient direction, $\partial L/\partial \mathbf{x}$, in steps of size $\tau$ to update the estimate. Using the objective function in Eq. (1), the update from step $k$ to $k$+1 can be expressed as follows,

$$\mathbf{x}_{k+1} = \mathbf{x}_k - \tau \mathbf{A}^T (\mathbf{A}\mathbf{x}_k - \mathbf{y}). \qquad (2)$$

Eq. (2) can be rearranged into the form of Richardson iteration [12] (Eq. (3)) by introducing a pre-calculated vector $\mathbf{b} = \tau \mathbf{A}^T \mathbf{y}$,

$$\mathbf{x}_{k+1} = (\mathbf{I} - \tau \mathbf{A}^T \mathbf{A})\mathbf{x}_k + \mathbf{b}, \qquad (3)$$

where each iteration computes two matrix-vector multiplications and two vector additions.

Due to the high-dimensional nature of practical inverse problems, repeated evaluation of Eq. (3) can be time-consuming. As an example, reconstructing a 256×256 image from full-rank measurements using Richardson iteration involves a matrix **A** with at least 65,536×65,536 elements. The resulting multiplication between **x** and **A** consists of more than 4 billion multiply-accumulate (MAC) operations. A conventional digital processor, such as a single CPU core, processes the matrix multiplication as sequential floating-point multiplication and addition operations, resulting in minutes to hours of computation time even for a moderate-scale imaging problem [13], [14].

An intuitive strategy for reducing the computation time is parallelizing the MAC operations, which is the fundamental principle behind numerous accelerators aiming to speed up the matrix multiplications. Digital accelerators, such as Tensor Core in the NVIDIA Volta GPU [15] and Google TPU [16], implement a large array of parallel, full-precision or half-precision floating-point processing units (FPUs) with application-specific integrated circuits (ASICs). Due to the on-chip power consumption limit, fixed-point digital multiply-adders that consume only a fraction of the energy per MAC than the FPUs [17], [18] become increasingly favorable in large-scale digital accelerators. Recently, analog accelerators that utilize natural physical processes for array multiplications and summations have shown to be even faster and more energy efficient for matrix processing [19], with implementations on both electronic [20], [21] and photonic [22], [23] platforms. However, analog accelerators are susceptible to signal artifacts arising from noises and device imperfections [24], and are unable to achieve floating-point precision. As a result, fixed-point is the de-facto data format in matrix accelerators [25]. Currently, these energy-efficient accelerators are limited to applications that are robust against computation errors, such as artificial neural networks [26]–[28]. For iterative solvers, floating-point precision is typically required [29], [30], because

This work was supported in part by the Army Research Office under Grant W911NF1910385, and in part by the Office of Naval Research under Grant N00014-20-1-2441.

Zheyuan Zhu, Andrew B. Klein, Guifang Li, and Shuo Pang are with CREOL, The College of Optics and Photonics, University of Central Florida, Orlando, FL 32816 USA (e-mails: zyzhu@knights.ucf.edu, andrewbklein@knights.ucf.edu, li@ucf.edu, pang@creol.ucf.edu).

low-precision matrix calculations can introduce cumulative error into the solution and stall the iterative process.

This work studies the properties of an iterative linear inverse algorithm in which the matrix computations are performed with fixed-point arithmetic. Using residual iterations for error correction, the precision of the fixed-point Richardson iteration can be theoretically extended to any desired level. Similar methods have been previously proposed for correcting the noise-induced errors in hybrid analog-digital processors [29], [30], but they lack quantitative analysis on the error and convergence, especially for fixed-point solvers. This paper is organized as follows. In Section 2, we introduce two theorems describing the convergence criteria and error of a fixed-point Richardson solver. The theorems have two major implications: 1) the fixed-point solver does not compromise the rate of convergence, and 2) using residual iterations, the fixed-point solver can reach solutions beyond its native precision limit. We verify the two theorems with three examples, namely, matrix inversion, image deconvolution, and tomography reconstruction in Section 3. Section 4 describes a prototype fixed-point iterative solver implemented on a field programable gate array (FPGA), which is more than two times faster than CPU but reaches the solution in the same number of iterations as a floating-point processor. Section 5 concludes the paper.

## II. THEORY

### A. Fixed-point format

In the fixed-point data format, a vector, matrix, or tensor is represented by an array of signed mantissas, while a global exponent is shared among all the elements. Each element, $x_i$, in the $L$-bit fixed-point representation, $\tilde{\mathbf{x}}$, of an $N$-element array $\mathbf{x} \in \mathbb{R}^N$ can be expressed as,

$$\tilde{x}_i = \text{sign}(x_i) \times \text{mant}_i \times 2^{\text{expo}-(L-1)}. \quad (4)$$

Here expo denotes the exponent indicating the number of bits above decimal point. To avoid overflow, the exponent is determined by the maximum element of the array,

$$\text{expo} = \left\lceil \log_2 \max_i(\text{abs}(x_i)) \right\rceil, \quad (5)$$

where $\lceil x \rceil$ denotes the rounding to the nearest integer larger than $x$. The $(L-1)$-bit mantissas $\text{mant}_i$ are calculated as,

$$\text{mant}_i = \lfloor abs(x_i) \times 2^{L-1-\text{expo}} \rfloor, \quad (6)$$

where $\lfloor x \rfloor$ denotes the rounding to the nearest integer smaller than $x$.

The multiplication between a fixed-point matrix $\widetilde{\mathbf{M}}$ and a fixed-point vector $\tilde{\mathbf{x}}$ can be simplified as integer arithmetic between the mantissas, accompanied by bit-shifting to match the exponent of the result. Compared with the floating-point matrix-vector product $\mathbf{y} = \mathbf{Mx}$, the $\ell^2$ norm of the error, $\mathbf{\delta y}$, of the fixed-point matrix-vector product $\tilde{\mathbf{y}} = \widetilde{\mathbf{M}}\tilde{\mathbf{x}}$ is bounded by,

$$|\mathbf{\delta y}|_2 = \left|\widetilde{\mathbf{M}}\tilde{\mathbf{x}} - \mathbf{Mx}\right|_2 \leq \eta |\mathbf{M}|_2 |\mathbf{x}|_2. \quad (7)$$

Here $|\mathbf{M}|_2$ denote the $\ell^2$ norm of the operator. An upper bound of the coefficient $\eta$ is given by $\zeta_v + \zeta_m + 3\zeta_v\zeta_m$ in Appendix I. The coefficients $\zeta_v$ and $\zeta_m$ are defined as:

$$\zeta_v := |\mathbf{x} - \tilde{\mathbf{x}}|_2/|\mathbf{x}|_2, \quad (8)$$
$$\zeta_m := \left|\mathbf{M} - \widetilde{\mathbf{M}}\right|_2/|\mathbf{M}|_2,$$

which represent the normalized rounding errors of the fixed-point vector and matrix, respectively. Both $\zeta_v$ and $\zeta_m$ can be determined from the bit widths and exponents of the arrays. In most cases, $\eta$ is much lower than this upper bound, depending on the distribution of elements in $\mathbf{M}$ and $\mathbf{x}$.

### B. Error and convergence criteria of fixed-point Richardson iteration

The solution $\mathbf{x}^*$ to a linear system $\mathbf{y} = \mathbf{Ax}$ is unique if matrix $\mathbf{A}$ has full rank. To ensure that the intermediate estimations in the Richardson iteration $\{\mathbf{x}_k\}$ converge to the solution $\mathbf{x}^*$ as the iteration number $k \to \infty$, all the eigenvalues of matrix $(\mathbf{I} - \tau\mathbf{A}^T\mathbf{A})$ must fall within $(-1,1)$. Since the largest and smallest eigenvalues are $1 - \tau|\mathbf{A}^T\mathbf{A}|_2/\kappa$ and $1 - \tau|\mathbf{A}^T\mathbf{A}|_2$, respectively, where the operator norm $|\mathbf{A}^T\mathbf{A}|_2$ equals the largest eigenvalue, and $\kappa$ is the condition number of $\mathbf{A}^T\mathbf{A}$, the choice of the step size, $\tau$, is confined by:

$$0 < \tau < \frac{2}{|\mathbf{A}^T\mathbf{A}|_2}. \quad (9)$$

In practice, to ensure the convergence while expediating the iterative solver, the maximum step size, $\tau_{max}$, is often selected with a safety margin $0 < \chi \ll 1$ as follows:

$$\tau_{max} = \frac{2-\chi}{|\mathbf{A}^T\mathbf{A}|_2}. \quad (10)$$

For fixed-point Richardson iterations, Theorem 1 below summarizes the error and convergence criteria based on the fixed-point matrix multiplication error $\eta$.

***Theorem 1***: For a linear system $\mathbf{Ax} = \mathbf{y}$ solved using fixed-point Richardson iterations (Eq. (3)) with step size $\tau$, under the conditions in Eq. (11),

$$0 < \tau < \frac{2}{|\mathbf{A}^T\mathbf{A}|_2} \text{ and } \eta < \frac{\tau|\mathbf{A}^T\mathbf{A}|_2}{\kappa - \tau|\mathbf{A}^T\mathbf{A}|_2}. \quad (11)$$

the asymptotic error, $\theta$, between the estimation, $\mathbf{x}_k$, and the solution, $\mathbf{x}^*$, is bounded by:

$$\theta = \lim_{k\to\infty} \sup \frac{|\mathbf{x}^* - \mathbf{x}_k|_2}{|\mathbf{x}^*|_2} \leq \eta\left(\frac{\kappa}{\tau|\mathbf{A}^T\mathbf{A}|_2} - 1\right), \quad (12)$$

where $\kappa$ is the condition number of $\mathbf{A}^T\mathbf{A}$, and $\eta$ is the normalized $\ell^2$ error of the two matrix-vector multiplications in each iteration.

Theorem 1 is proved in Appendix II by examining the $\ell^2$ norm of the error of intermediate estimations, $|\mathbf{x}^* - \mathbf{x}_k|_2$, as $k \to \infty$. The proof shows a linear convergence with rate $|\mathbf{I} - \tau\mathbf{A}^T\mathbf{A}/\kappa|_2$, which is independent on the error of the fixed-point matrix-vector product, $\eta$. Because a typical problem involves a value of $\kappa$ on the order of $10^2 \sim 10^3$, using the choice of step size $\tau$ in Eq. (10), $\tau|\mathbf{A}^T\mathbf{A}|_2, = 2-\chi \ll \kappa$, we can approximate

$$|\mathbf{I} - \tau\mathbf{A}^T\mathbf{A}/\kappa|_2^k \approx \exp(-\gamma k), \quad (13)$$

where the convergence rate $\gamma = (2-\chi)/\kappa$.

### C. Fixed-point residual iteration

The asymptotic error $\theta$ in Theorem 1 is independent of the solution $\mathbf{x}^*$, which implies that the fixed-point iterative solver can stagnate at an estimation, $\mathbf{x}_\infty$, closer to the solution $\mathbf{x}^*$ in terms of $\ell^2$ distance, if $\mathbf{x}^*$ has a lower $\ell^2$ norm. Hence, constructing a sequence of fixed-point iterative solvers in which the solutions have decreasing $\ell^2$ norm becomes the natural choice to approach the unbiased solution of the inverse problem.

ALGORITHM I
FIXED-POINT RESIDUAL ITERATION FOR SOLVING THE LINEAR SYSTEM $\mathbf{y} = \mathbf{A}\mathbf{x}$

| Residual iteration algorithm |
|---|
| $\mathbf{x}^{(0)} \leftarrow 0,\ \mathbf{r}^{(0)} \leftarrow \mathbf{y},$ <br> for $l$ from 1 to $M$ <br>     $\delta \mathbf{x}_0 \leftarrow \mathbf{0}$ <br>     $\mathbf{b}^{(l)} \leftarrow \tau \mathbf{A}^T \mathbf{r}^{(l)}$ <br>     // Solve $\mathbf{r}^{(l)} = \mathbf{A}\delta\mathbf{x}^{(l)}$ with fixed-point Richardson iterations <br>     repeat: <br>         $\delta\mathbf{x}_k \leftarrow \delta\mathbf{x}_{k-1} - \mathbf{A}^T(\mathbf{A}\delta\mathbf{x}_{k-1}) + \mathbf{b}^{(l)}$ <br>     until $|\delta\mathbf{x}_k - \delta\mathbf{x}_{k-1}|_2 < \varepsilon.$ <br>     $\delta\mathbf{x}^{(l)} \leftarrow \delta\mathbf{x}_k$ <br>     // Solution and residue update <br>     $\mathbf{x}^{(l)} \leftarrow \mathbf{x}^{(l-1)} + \delta\mathbf{x}^{(l)}$ <br>     $\mathbf{r}^{(l)} \leftarrow \mathbf{y} - \mathbf{A}\mathbf{x}^{(l)}$ |

The proposed residual iteration algorithm described in Algorithm 1 leverages the adjustable exponent of the fixed-point format. The algorithm contains $M$ sets of fixed-point Richardson iterations as the inner loops. The acceptance criterion of the inner loop solution is:
$$|\delta\mathbf{x}_k - \delta\mathbf{x}_{k-1}|_2 < \varepsilon, \tag{14}$$
where $\varepsilon$ is a small number of the order $\zeta_v$. Let $\mathbf{x}^{(M)} = \sum_{l=1}^{M} \delta\mathbf{x}^{(l)}$ denote the accumulated solutions from $M$ fixed-point Richardson iterative solvers. The proposed algorithm aims to reduce the norm of $\mathbf{r}^{(l)} \coloneqq \mathbf{y} - \mathbf{A}\mathbf{x}^{(l)}$ at each outer loop $l$ by recursively solving the linear system $\mathbf{r}^{(l)} = \mathbf{A}\delta\mathbf{x}^{(l)}$ with the updated exponents in the next inner loop of the fixed-point Richardson iterations.

It is worth noting that similar precision refinement methods [31], [32] have been employed in iterative solvers using mixed full- and half-precision floating-point data formats. Yet, the floating-point precision refinement typically requires far fewer residual loops, $M$, than Algorithm 1. In addition, floating-point format allows different exponents for every element in the matrix-vector product. As a result, the error analysis in floating-point precision refinement [33] does not apply to fixed-point residual iterations. Here we present the upper bound of the normalized error of the solution $\mathbf{x}^{(M)}$ after $M$ fixed-point residual iterations in Theorem 2.

***Theorem 2***: For a full-rank linear system $\mathbf{A}\mathbf{x} = \mathbf{y}$ with its solution $\mathbf{x}^*$ obtained from the residual iteration algorithm (Algorithm 1), the asymptotic error of the estimation after $M$ residue updates, $\mathbf{x}^{(M)}$, is bounded by:
$$\theta^{(M)} \leq \theta^M, \tag{15}$$
where $\theta$ is the asymptotic error of the solution obtained from fixed-point Richardson solver (Eq. (12)).

Theorem 2 is established based on the decreasing $\ell^2$ norm of the solution $\mathbf{x}^{*(l)}$ in each fixed-point Richardson iteration. The upper bound of the solution error after $M$ residue updates can be derived recursively from Eq. (12) in Theorem 1. Detailed proof is presented in Appendix III.

For Eq. (15) to be valid, the exponents of $\mathbf{b}^{(l)}$ and $\mathbf{A}^T\mathbf{A}\delta\mathbf{x}_k$ should be updated after every iteration of the inner loop. In actual implementations shown in Section 3, due to the computational cost of distribution-based estimation of the infinity norm (Eq. (19)), the exponents are adjusted after the completion of each inner loop. As a result, the errors may not always satisfy the bound in Eq.(15).

### III. FIXED-POINT ITERATIVE SOLVER EXPERIMENTS

#### A. Matrix inversion

We first verify the convergence rate and the residue error in theorem 1 with a fixed-point matrix inversion solver. The two matrices to be inverted, $\mathbf{A}_1$ and $\mathbf{A}_2$, shown in Fig. 1 (a), are constructed from the 4×4 discrete cosine transfer matrix by applying different linear scaling factors to its eigenvalues. The condition numbers $\kappa$ of $\mathbf{A}_1$ and $\mathbf{A}_2$ are 25.0 and 11.1, respectively. The analytical inversions $\mathbf{A}_1^{-1}$ and $\mathbf{A}_2^{-1}$, shown in Fig. 1 (b), are calculated with floating-point LU decomposition.

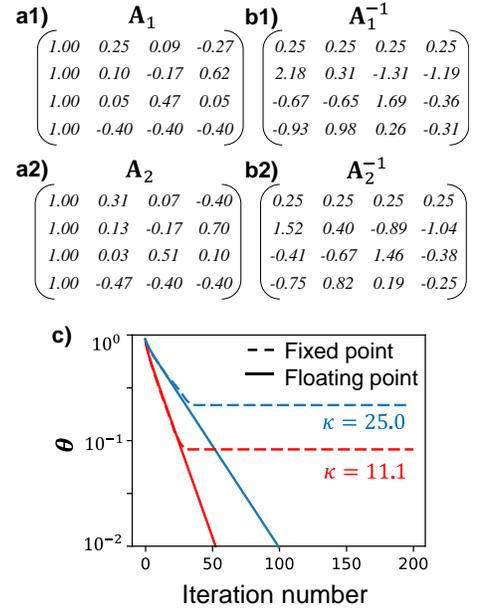

Fig. 1. Two 4×4 matrices (a) and their inverses (b) in fixed-point iterative matrix inversion solver. $\mathbf{A}_1$ has a $\kappa$ of 25, and $\mathbf{A}_2$ has a $\kappa$ of 11.1. All matrices are stored in floating-point format but displayed with two decimal points. (c) Normalized error of the iterative solutions $\theta_k$ vs. the iteration number $k$ for the inversion of $\mathbf{A}_1$ ($\kappa$=25.0) and $\mathbf{A}_2$ ($\kappa$=11.1) using 8-bit fixed-point Richardson iterations.

TABLE 1. NORMALIZED ERRORS AND THE CONVERGENCE RATE OF THE 8-BIT RICHARDSON ITERATIONS USED TO OBTAIN THE INVERSE OF $\mathbf{A}_1$ AND $\mathbf{A}_2$. A 95% CONFIDENCE BOUND IS USED IN EXPONENTIAL FITTING.

| Matrix to invert | $\kappa$ | $\theta$ | Upper bound of $\theta$ | Convergence rate $\gamma$ |
|---|---|---|---|---|
| $\mathbf{A}_1$ | 25.0 | 0.21 | 0.24 | 0.040±0.001 |
| $\mathbf{A}_2$ | 11.1 | 0.083 | 0.098 | 0.085±0.004 |

The solver uses signed 8-bit fixed-point format. The step sizes, $\tau$, are calculated from Eq. (10) with a safety margin, $\chi$=0.2, for the inversions of $\mathbf{A}_1$ and $\mathbf{A}_2$. To estimate the normalized $\ell^2$ error $\eta$ of the matrix-vector product, $\mathbf{A}^T\mathbf{A}\mathbf{x}_k$, the arrays $\mathbf{A}^T\mathbf{A}$ and $\mathbf{b}$ in Eq. (3) are pre-calculated and converted to fixed-point. We enumerate 200 normally distributed $\mathbf{x}$ with mean $\mathbf{b}$ and standard deviation $|\mathbf{b} - \mathbf{x}^*|_2/2$, where $\mathbf{x}^*$ is the solution. The enumerated $\mathbf{x}$ are converted to fixed-point format and multiplied with the fixed-point matrix $\widetilde{\mathbf{A}}^T\widetilde{\mathbf{A}}$. We then obtain

$\eta$=0.019 from the average normalized root-mean-square error (RMSE) of the multiplication results.

Fig. 1 (c) plots the normalized errors, $\theta_k$, in log-scale compared with the analytical inverses, $\mathbf{A}_1^{-1}$ and $\mathbf{A}_2^{-1}$, as the iterations progress. The asymptotic errors, $\theta$, of the 8-bit iterations are 0.21 and 0.083 for the inversions of $\mathbf{A}_1$ and $\mathbf{A}_2$, respectively, both of which are consistent with the error bounds given by Eq. (12). The convergence rate, $\gamma$, is obtained from the exponential fitting of Fig. 1 (c) with a 95% confidence bound. For the inversions of $\mathbf{A}_1$ and $\mathbf{A}_2$, $\gamma$ is 0.040±0.001 and 0.085±0.004 respectively, which are inversely proportional to $\kappa$, as predicted by Eq. (13).

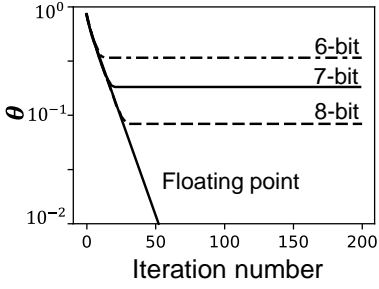

Fig. 2. Normalized errors of the iterative solutions $\theta_k$ vs. the iteration number $k$ for the inversion of $\mathbf{A}_2$ ($\kappa$=11.1) using 8-bit, 7-bit, and 6-bit fixed-point iterative solvers. The error of the floating-point iterative solver is plotted for reference.

TABLE 2. ASYMPTOTIC ERRORS AND CONVERGENCE RATES OF THE 8-BIT, 7-BIT, 6-BIT, AND FLOATING-POINT ITERATIONS. A 95% CONFIDENCE BOUND IS USED IN EXPONENTIAL FITTING.

| Precision | $\eta$ | $\theta$ | Upper bound of $\theta$ | Convergence rate $\gamma$ |
|---|---|---|---|---|
| 8-bit | 0.019 | 0.083 | 0.098 | 0.085±0.002 |
| 7-bit | 0.036 | 0.18 | 0.19 | 0.085±0.004 |
| 6-bit | 0.072 | 0.33 | 0.37 | 0.088±0.009 |
| Floating-point | N.A. | N.A. | N.A. | 0.083±0.003 |

For the same condition number, Fig. 2 compares the fixed-point Richardson iterations with the 8-bit, 7-bit, 6-bit, and floating-point precisions. The convergence rates are summarized in Table 2, which confirms the implication of Eq. (13) that the fixed-point iterative solver maintains the same convergence rate, $\gamma$, as the floating-point counterpart. The independence of convergence rate on the precision suggests that fixed-point iterative solvers do not compromise the convergence speed. The asymptotic errors of the 8-bit, 7-bit, and 6-bit solutions are 0.083, 0.18, and 0.33, respectively, which are all below the upper bound established in Eq. (12).

*B. Image deconvolution*

Deconvolution is a classic linear inverse problem that has broad applications in signal/image filtering and processing. For a linear, shift-invariant system, the convolution output, $\mathbf{y}$ (i.e., filtered signal, or blurred image) is the discrete convolution between the input signal or image $\mathbf{x}$ and the kernel $\mathbf{K}$, i.e., $\mathbf{y} = \mathbf{K} * \mathbf{x}$, which is defined elementwise as:

$$y_{i,j} = \sum_{p,q=0}^{M,N} K_{p,q} x_{i-p,j-q}, \qquad (16)$$

where $M$ and $N$ are the sizes of the kernel along row and column directions, respectively. Here, we have assumed a two-dimensional convolution, as in the case of image processing. This example demonstrates the iterative Richardson-Lucy deconvolution, which recovers images from blurry measurements given the convolution kernel.

The discrete convolution can be formulated as the product between the vectorized input, $\mathbf{x}$, and a Toeplitz matrix, $\mathbf{A}$, constructed from the kernel [34], producing the output $\mathbf{y} = \mathbf{A}\mathbf{x}$ in the vectorized form. Using the convolution theorem, Eq. (16) can be rewritten as

$$\mathbf{y} = \mathcal{F}^\dagger \mathrm{diag}(\widetilde{\mathbf{K}}) \mathcal{F} \mathbf{x}, \qquad (17)$$

where $\mathcal{F}$ and $\mathcal{F}^\dagger$ represent the forward and inverse two-dimensional Fourier transforms, both of which are unitary operators, † denotes the conjugate transpose of a complex matrix, and $\mathrm{diag}(\widetilde{\mathbf{K}})$ is the diagonal matrix constructed from the Fourier transform of the kernel $\widetilde{\mathbf{K}}$. The condition number, $\kappa$, of $\mathbf{A}^T\mathbf{A}$ can be obtained from the maximum and minimum elements in $\mathrm{diag}(\widetilde{\mathbf{K}})^\dagger \mathrm{diag}(\widetilde{\mathbf{K}})$.

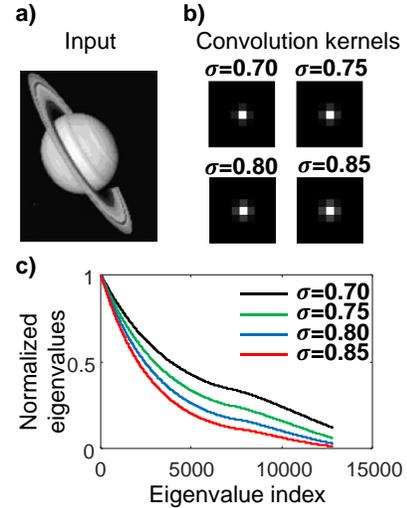

Fig. 3. Fixed-point iterative solver for a discrete Richardson-Lucy deconvolution problem. (a) The 128×102 signed 4-bit input image "Saturn". (b) Four convolution kernels with $\sigma$=0.70, 0.75, 0.80, and 0.85, respectively. (c) Eigenvalue spectra of the Toeplitz matrices calculated from the Fourier transforms of the convolution kernels. The condition numbers, $\kappa$, of $\mathbf{A}^T\mathbf{A}$ are 8.38, 16.59, 35.07, and 78.50, respectively.

Fig. 3 shows an example based on Richardson-Lucy deconvolution. The input, $\mathbf{x}^*$, is a test image downsampled from the MATLAB image "Saturn", shown in Fig. 3 (a). The image contains 128×102 pixels in signed 4-bit fixed-point format with range (-1,1). The convolution kernel, $\mathbf{K}$, follows the Gaussian profile,

$$K_{p,q} = K_0 \exp\left(-\frac{\left(p - \frac{M}{2}\right)^2 + \left(q - \frac{N}{2}\right)^2}{\sigma^2}\right), \qquad (18)$$

where the parameter $\sigma$ determines the width of the kernel, and $K_0$ is a normalization constant that ensures that $\sum_{p,q}|K_{p,q}|^2 =$

1. Four example kernels with $\sigma$=0.70, 0.75, 0.80, and 0.85, are plotted in Fig. 3 (b). Fig. 3 (c) shows the normalized eigenvalue spectra of the four Toeplitz matrices, calculated from the Fourier transform of the kernels. The condition numbers, $\kappa$, of the four kernels are 8.38, 16.59, 35.07, and 78.50, respectively. If $\sigma$ is increased beyond 0.90, $\kappa$ would violate the criteria in Eq. (11) and lead to a non-converging deconvolution result. Similar phenomenon has been reported in hybrid analog-digital solvers [35].

The measurements, shown in Fig. 4 (a), are calculated from convolving the input image, $\mathbf{x}^*$, with the four kernels in Fig. 3 (b), and then quantized to 8-bit to simulate the digitization error from a CCD or CMOS camera. The deconvolutions use 8-bit fixed-point iterations, with $\eta$ estimated to be 0.014. Step sizes with the safety margin $\chi$=0.2 are used for all four iterative deconvolutions. Fig. 4 (b) plots the deconvolved images after 200 iterations. The asymptotic errors, $\theta$, are summarized in Table 3, and are all below the error bound predicted in Eq. (12).

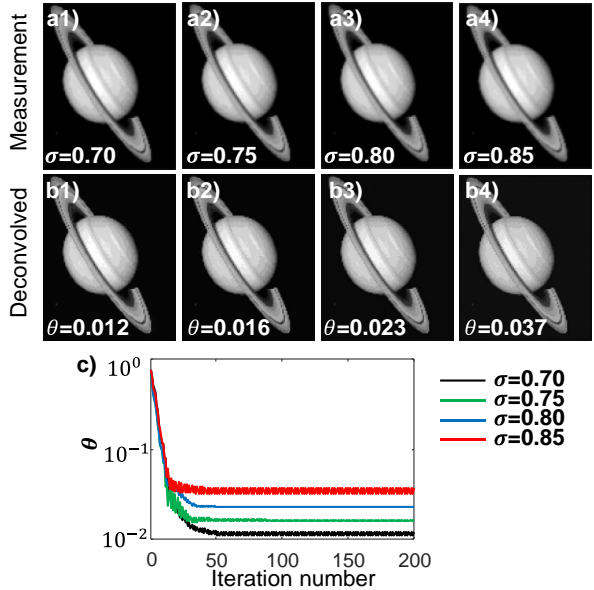

Fig. 4. 8-bit fixed-point Richardson solver of a discrete Richardson-Lucy deconvolution problem. (a) Convolution measurements from the four kernels with $\sigma$ of (a1) 0.70, (a2) 0.75, (a3) 0.80, and (a4) 0.85. (b1-b4) Deconvolved images using 8-bit fixed-point iterations from (a1-a4). (c) Normalized error vs. iteration steps for the reconstructions in (b1-b4).

TABLE 3. NORMALIZED ERRORS OF RICHARDSON-LUCY DECONVOLUTION RESULTS.

| $\sigma$ | $\kappa$ | $\theta$ | Upper bound of $\theta$ |
|---|---|---|---|
| 0.70 | 8.38 | 0.012 | 0.056 |
| 0.75 | 16.59 | 0.016 | 0.12 |
| 0.80 | 35.07 | 0.023 | 0.25 |
| 0.85 | 78.50 | 0.037 | 0.53 |

C. *Tomographic reconstruction with residual iteration*

Theorem 2 indicates that the residual iteration can overcome the stagnation of the fixed-point Richardson solver, and reach a solution beyond the fixed-point limit. To verify the performance of the residual iterations, we examine a tomographic reconstruction of a 16×16 "Super Mario" pixel art from its pencil-beam projections, shown in Fig. 5 (a). The pixels of the input image, $\mathbf{x}^*$, are stored in the signed 3-bit fixed-point format with range (-1,1), giving a quantization interval of 0.25 between two adjacent intensity levels. The measurements consist of 60 projections between 0 and 180º at 4º intervals, each containing 31 uniformly spaced beams. The system matrix, $\mathbf{A}$, of the tomographic model is constructed from the distance-driven method [36]. Fig. 5 (b) plots the eigenvalue spectrum of the tomographic model. The condition number, $\kappa$, of $\mathbf{A}^T\mathbf{A}$ is 101.80. Assuming that the maximum step size, $\tau$, is used in the Richardson iterations, the computation error $\eta$ must be less than 0.020 to ensure convergence, according to Eq. (11). For the size and sparsity of this system matrix, we find that a minimum of $L$=8 bits (with an estimated $\eta$=0.015) are required for fixed-point iterations to converge.

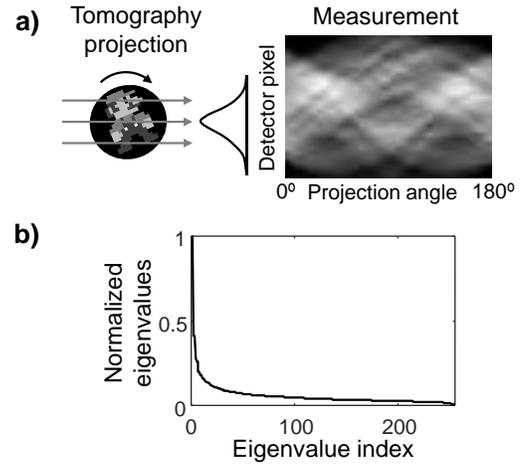

Fig. 5. Fixed-point iterative solver for a tomographic reconstruction problem. (a) Tomographic projection and measurement of a 16×16 signed 3-bit "Super Mario" pixel art. (b) Eigenvalue spectrum of the tomographic projection model. The condition number $\kappa$ of $\mathbf{A}^T\mathbf{A}$ is 101.80.

Fixed-point Richardson and residual iterations with 8-bit, 9-bit, and 10-bit precisions are performed to reconstruct the "Super Mario" pixel art from the tomographic projections. The step sizes, $\tau$, with safety margin $\chi$=0.3 are used in all the solvers. The system matrix $\mathbf{A}$ is quantized to the solver precision with expo=0. The measurement, $\mathbf{y}$, is serialized to a vector, multiplied by $\tau\mathbf{A}^T$, and quantized to the solver precision with expo=0 to pre-calculate the fixed-point vector $\mathbf{b}$ for the iterations. The fixed-point Richardson iteration sets the exponent of the multiplication between $\mathbf{A}^T\mathbf{A}$ and $\mathbf{x}_k$ to 2. For residual iterations, the exponents of matrix-vector multiplications are decremented by 1 every 80 steps from 2 to -2.

Fig. 6 plots the reconstructions from Richardson and residual iterations with 8-bit, 9-bit, and 10-bit precisions. Table 4 summarizes the normalized errors of the fixed-point Richardson reconstructions, all of which are below the error bound given by Eq. (12). The normalized error, $\theta$, is plotted as a function of the Richardson iteration number in Fig. 6 (c1). Note that the same convergence rate is present regardless of the precision of the solver. The normalized errors of the residual iterations are shown in Fig. 6 (c2). After $M$=5 residue updates, the 8-bit, 9-

bit, and 10-bit solvers all reach normalized errors below 0.1. At this level of $\theta$, the pixelwise differences between $\mathbf{x}^{(5)}$ and $\mathbf{x}^*$ are below the quantization interval, 0.25, of the true image $\mathbf{x}^*$. Therefore, the reconstructions in Fig. 6 (b1-b3) would exhibit no visual differences from the true image if displayed in 3-bit grayscale colormap.

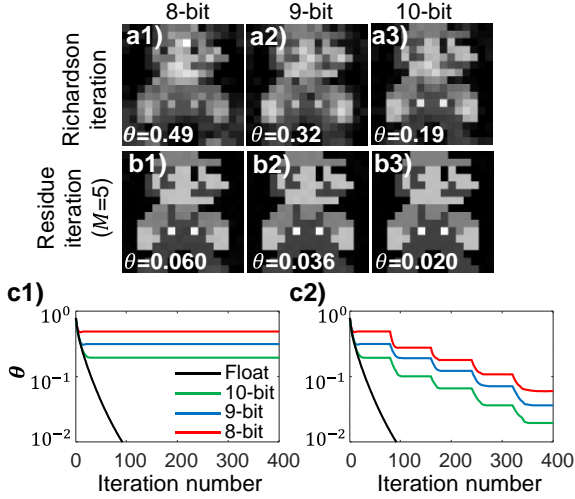

Fig. 6. Richardson (a) and residual (b) iterations for reconstructing the "Super Mario" pixel art from tomographic measurements. (a) Reconstructions from 8-bit (a1), 9-bit (a2), and 10-bit (a3) Richardson iterations. (b) Reconstructions from 8-bit (b1), 9-bit (b2), and 10-bit (b3) residual iterations after $M=5$ residue update steps. (c) Normalized error vs. iteration steps for 8-bit (red curves), 9-bit (blue curves), and 10-bit (green curves) Richardson (c1) and residual (c2) iterations. The floating-point (black curve) reconstruction is plotted on both (c1) and (c2) for reference.

TABLE 4. NORMALIZED ERRORS OF 8-BIT, 9-BIT, AND 10-BIT RICHARDSON ITERATIONS.

| Precision | $\eta$ | $\theta$ | Upper bound of $\theta$ |
|---|---|---|---|
| 8-bit | 0.015 | 0.49 | 0.88 |
| 9-bit | 0.0056 | 0.32 | 0.33 |
| 10-bit | 0.0028 | 0.19 | 0.16 |

Faster convergence is possible via adaptive adjustments of the exponents of $\mathbf{b}^{(l)}$ and $\mathbf{A}^T\mathbf{A}\delta\mathbf{x}_k$. Fig. 7 shows the error of the 8-bit residual iteration executed on the prototype solver detailed in section 4. The exponent is updated every 5 steps based on the distribution of the elements in a fixed-point array $\mathbf{x}$. The exponents are calculated using:

$$\text{expo} = \lceil \log_2(|\mu_\mathbf{x}| + 3\sigma_\mathbf{x}) \rceil, \tag{19}$$

where $\mu_\mathbf{x}$ and $\sigma_\mathbf{x}$ denote the mean and standard deviation, respectively, of all the elements in $\mathbf{x}$. The adaptive residual iteration achieves the same convergence speed as the floating-point solver, with both methods achieving errors below the quantization interval of the true solution after 28 iterations.

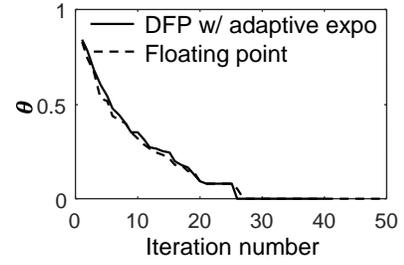

Fig. 7. Normalized error vs. iteration number for the floating-point solver (dashed curve) and fixed-point iterative solver (solid curve), with adaptive exponent adjustments made every 5 steps.

## IV. HARDWARE PROTOTYPE FOR A FIXED-POINT ITERATIVE SOLVER

We have built a hardware prototype to perform the fixed-point Richardson iteration and exponent adjustment. The core component of the solver, i.e., the fixed-point matrix multiplication, is implemented on an FPGA development board (Xilinx ZC706) and communicates with the host PC via PCIe bus. Fig. 8 illustrates the block diagram of the logic design within the FPGA. The systolic multiplier array performs 512 MAC operations in a single clock cycle. The inputs of the multiplier array consist of a 16×16 block $\mathcal{W}$ and a 16×2 block $\mathcal{X}$. The precisions of both input blocks are signed 8-bit. The multiplication outputs are cached in signed 32-bit and wired to a 16×2 array of accumulators. For the matrix-vector multiplication involving a matrix $\mathbf{W}$ larger than 16×16, the matrix $\mathbf{W}$ is zero-padded to an integer number of 16×16 blocks $\mathcal{W}$, and the batch vector $\mathbf{X}$ is zero-padded to an integer number of 16×2 blocks $\mathcal{X}$, along the row dimension only. All the decomposed blocks $\mathcal{W}$ and $\mathcal{X}$ are store in cache. The multiplier array reads the $\mathcal{W}$ blocks along the rows and the $\mathcal{X}$ blocks along the columns. The blocks of 16×2 results are summed in the accumulators. After completion of one row of $\mathcal{W}$ blocks, the multiplier array moves to the next row block of $\mathcal{W}$ and cycles through all the blocks of $\mathcal{X}$ again.

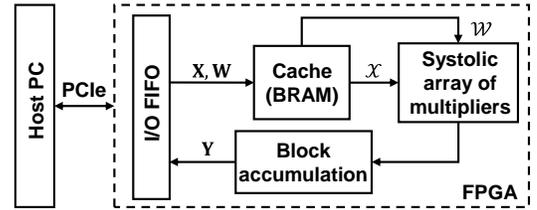

Fig. 8. Block diagram of the fixed-point iterative solver prototype on FPGA depicting the communications and among systolic multiplier array, cache, accumulator, and other resources. FIFO: first-in-first-out; BRAM: block random-access memory.

The host PC converts $\mathbf{W}$ and $\mathbf{X}$ into fixed-point, and then decomposes the mantissas into 16×16 $\mathcal{W}$ blocks and 16×2 $\mathcal{X}$ blocks. The decomposed $\mathcal{W}$ and $\mathcal{X}$ blocks are streamed to the FPGA through PCIe. Once all the blocks are multiplied and accumulated in FPGA, the results $\mathbf{Y}$ are streamed back to host PC in 32-bit integer format. The host PC masks the results back to 8-bit, with the most significant bit (MSB) and least significant bit (LSB) selected by the exponents. Every five iterations, the host PC updates the exponents according to Eq.

(19).

For the FPGA design in Fig. 8 implemented on the ZC706 evaluation board, the Xilinx synthesis and implementation tool estimates a power consumption of 2.047W, of which 1.253W is consumed by the transceivers for PCIe communication with the host PC, and the remaining 0.794W is consumed by the logics, including the clocks, systolic multiplier array, BRAM cache, and control units that generate the addresses of $\mathcal{W}$ and $\mathcal{X}$ blocks. The clock rate is 250MHz. Considering the total number of 512 MACs within a single clock cycle, the energy consumption by the fixed-point computation is 6.2pJ/MAC, which is on the same order as that of in-datacenter TPUs [37].

We have tested the speed of the iterative solver implemented on our hardware prototype and CPU (Intel Core i3-4130, dual-core, 3.40GHz clock rate). The total CPU time of the fixed-point Richardson iterations, excluding exponent adjustment, is measured, and then divided by the total number of iterations in all residual updates to obtain the time per iteration, which is indicative of the calculation speed. The fixed-point solver on CPU takes 1.7±0.2ms per iteration, while that on FPGA prototype takes 0.76±0.05ms per iteration, representing more than two times of improvement in speed. Errors represent the variance of the calculation time for five repeated tests.

## V. CONCLUSION

We have demonstrated a fixed-point iterative solver that computes high-precision solutions for linear inverse problems beyond the precision limit of the hardware. We have described the convergence and error of the fixed-point iterations in two theorems. Theorem 1 shows that the convergence rate is independent of the precision of the solver, which implies that fixed-point iteration does not compromise the convergence speed of the solver. Theorem 2 resolves the stagnation of a fixed-point solver with residual iterations, which correct the error and refine the solution to a higher precision level. We have presented three examples of linear inverse problems, namely, matrix inversion, Richardson-Lucy deconvolution, and tomographic reconstruction, to verify both theorems.

The combination of residual iteration with adaptive exponent adjustment achieves the same rate constant and error as a floating-point Richardson solver. We have prototyped the fixed-point solver with a fixed-point systolic multiplier array on the FPGA and fixed-point arithmetic libraries on host PC. The prototype solver achieves more than twice the calculation speed of a CPU and takes the same number of iterations to reach the solution as its floating-point counterpart. We envision that a low-precision, analog matrix processing core can supersede the current FPGA systolic array in the next generation of fast, energy-efficient solvers for numerous large-scale, high-performance computing applications, including but not limited to convex optimizations, differential equations, and training of artificial neural networks.

## APPENDIX I

We derive the upper bound of the normalized error, $\eta$ for fixed-point matrix-vector multiplication. Based on Eq. (5), for an $N$-dimensional vector, $\mathbf{x}$, the infinity norm falls in the range:
$$2^{\text{expo}-1} \le |\mathbf{x}|_\infty \le 2^{\text{expo}}. \tag{20}$$
The infinity norm of fixed-point rounding error, $|\mathbf{x} - \tilde{\mathbf{x}}|_\infty$, is determined by the bit width:
$$0 \le |\mathbf{x} - \tilde{\mathbf{x}}|_\infty \le 2^{\text{expo}-(L-1)}. \tag{21}$$
The equivalence of norms establishes the inequality between infinity norm and $\ell^2$ norm, $|\mathbf{x}|_\infty \le |\mathbf{x}|_2 \le \sqrt{N}|\mathbf{x}|_\infty$. This gives the normalized $\ell^2$ error, $\zeta_v = |\mathbf{x} - \tilde{\mathbf{x}}|_2/|\mathbf{x}|_2$, of a fixed-point vector $\tilde{\mathbf{x}}$,
$$0 \le \zeta_v \le 2^{2-L}\sqrt{N}. \tag{22}$$
For an $N \times N$ matrix $\mathbf{M}$, its operator norm is used in error analysis. Using the same fixed-point conversions (Eq. (4)-(6)), the infinity norm of a fixed-point matrix $|\mathbf{M}|_\infty$ falls in the range:
$$2^{\text{expo}-1}N \le |\mathbf{M}|_\infty \le 2^{\text{expo}}N. \tag{23}$$
The infinity norm of the rounding error, $\left|\mathbf{M} - \widetilde{\mathbf{M}}\right|_\infty$, is determined by both bit width and size, and falls in the range:
$$0 \le \left|\mathbf{M} - \widetilde{\mathbf{M}}\right|_\infty \le 2^{\text{expo}-(L-1)}N. \tag{24}$$
The equivalence of norms for matrix yields $|\mathbf{M}|_\infty/\sqrt{N} \le |\mathbf{M}|_2 \le \sqrt{N}|\mathbf{M}|_\infty$ [29], providing the boundaries of the normalized $\ell^2$ error, $\zeta_m = \left|\mathbf{M} - \widetilde{\mathbf{M}}\right|_2/|\mathbf{M}|_2$, of a fixed-point matrix $\widetilde{\mathbf{M}}$:
$$0 \le \zeta_m \le 2^{2-L}N. \tag{25}$$
We introduce the rounding errors of the fixed-point matrix $\mathbf{E} = \mathbf{M} - \widetilde{\mathbf{M}}$ and the vector $\mathbf{e} = \mathbf{x} - \tilde{\mathbf{x}}$. Using Eq. (22) and Eq. (25), the $\ell^2$ norm of the error $\delta\mathbf{y} = \widetilde{\mathbf{M}}\tilde{\mathbf{x}} - \mathbf{M}\mathbf{x}$ is bounded by,
$$\begin{aligned}|\delta\mathbf{y}|_2 &\le \left|\widetilde{\mathbf{M}}\mathbf{e}\right|_2 + |\mathbf{E}\tilde{\mathbf{x}}|_2 + |\mathbf{E}\mathbf{e}|_2 \\ &\le \left|\widetilde{\mathbf{M}}\right|_2|\mathbf{e}|_2 + |\mathbf{E}|_2|\tilde{\mathbf{x}}|_2 + |\mathbf{E}|_2|\mathbf{e}|_2 \\ &= \zeta_v\left|\widetilde{\mathbf{M}}\right|_2|\mathbf{x}|_2 + \zeta_m|\mathbf{M}|_2|\tilde{\mathbf{x}}|_2 + \zeta_v\zeta_m|\mathbf{M}|_2|\mathbf{x}|_2 \\ &\le \zeta_v(1+\zeta_m)|\mathbf{M}|_2|\mathbf{x}|_2 + \zeta_m(1+\zeta_v)|\mathbf{M}|_2|\mathbf{x}|_2 \\ &\quad + \zeta_v\zeta_m|\mathbf{M}|_2|\mathbf{x}|_2 \\ &= (\zeta_v + \zeta_m + 3\zeta_v\zeta_m)|\mathbf{M}|_2|\mathbf{x}|_2.\end{aligned} \tag{26}$$
Let us define the normalized error of the matrix-vector product as,
$$\eta \coloneqq \frac{|\delta\mathbf{y}|_2}{|\mathbf{M}|_2|\mathbf{x}|_2}, \tag{27}$$
the upper bound of $\eta$ is $\zeta_v + \zeta_m + 3\zeta_v\zeta_m$.

## APPENDIX II

We prove Theorem 1 by examining the $\ell^2$ error of the intermediate estimates $|\mathbf{x}^* - \mathbf{x}_k|_2$. Without loss of generality, the normalized $\ell^2$ error $\eta$ of two matrix-vector multiplications can be combined into one multiplication error $\delta\mathbf{y}$ between a matrix, $\mathbf{B} \coloneqq \mathbf{I} - \tau\mathbf{A}^T\mathbf{A}$, and the vector $\mathbf{x}_k$. Since $\mathbf{x}^* = \mathbf{B}\mathbf{x}^* + \mathbf{b}$, the error of the solution at iteration $k$ can be obtained recurrently,
$$\begin{aligned}\mathbf{x}^* - \mathbf{x}_k &= \mathbf{B}(\mathbf{x}^* - \mathbf{x}_{k-1}) - \delta\mathbf{y} \\ &= \mathbf{B}^2(\mathbf{x}^* - \mathbf{x}_{k-2}) - (\mathbf{B}\delta\mathbf{y} + \delta\mathbf{y}) \\ &= \cdots \\ &= \mathbf{B}^k(\mathbf{x}^* - \mathbf{x}_0) - \sum_{p=0}^{k-1}(\mathbf{B})^p\delta\mathbf{y}.\end{aligned} \tag{28}$$
Assuming the initial guess, $\mathbf{x}_0 = \mathbf{0}$, the $\ell^2$ norm of the error at iteration $k$ is,

$$|\mathbf{x}^* - \mathbf{x}_k|_2 \leq |\mathbf{B}^k \mathbf{x}^*|_2 + \left|\sum_{p=0}^{k-1}(\mathbf{B})^p \delta \mathbf{y}\right|_2$$

$$\leq |\mathbf{B}|_2^k |\mathbf{x}^*|_2 + \sum_{p=0}^{k-1}|(\mathbf{B})|_2^p |\delta \mathbf{y}|_2 \quad (29)$$

$$\leq |\mathbf{B}|_2^k |\mathbf{x}^*|_2 + \frac{1-|\mathbf{B}|_2^k}{1-|\mathbf{B}|_2}|\delta \mathbf{y}|_2.$$

Substituting the matrix-vector product error $|\delta\mathbf{y}|$ with the upper bound $\eta|\mathbf{B}||\mathbf{x}^*|$ as in Eq. (7), the normalized error of the solution, $\theta_k$, at iteration $k$ is bounded by,

$$\theta_k = \frac{|\mathbf{x}^* - \mathbf{x}_k|_2}{|\mathbf{x}^*|_2}$$
$$\leq \frac{\eta|\mathbf{B}|_2}{1-|\mathbf{B}|_2} + \left(1 - \frac{\eta|\mathbf{B}|_2}{1-|\mathbf{B}|_2}\right)|\mathbf{B}|_2^k, \quad (30)$$

where the first term is the asymptotic error, and the second term is a function of the iteration step $k$. To ensure $\theta_k$ decays as the iteration progresses, the convergence criteria of the Richardson iteration must be satisfied,

$$\begin{cases} -1 < |\mathbf{B}|_2 < 1 \\ 1 - \frac{\eta|\mathbf{B}|_2}{1-|\mathbf{B}|_2} > 0 \end{cases}. \quad (31)$$

The largest and the smallest eigenvalues in $\mathbf{B}$ are $1-\tau|\mathbf{A}^T\mathbf{A}|_2/\kappa$ and $1-\tau|\mathbf{A}^T\mathbf{A}|_2$, respectively. Plugging $|\mathbf{B}|_2$ in the convergence criteria (Eq.(31)) confirms Eq. (11). Using the maximum eigenvalue $|\mathbf{B}|_2 = 1-\tau|\mathbf{A}^T\mathbf{A}|_2/\kappa$ in the first term of Eq. (30) gives the upper bound of asymptotic error $\theta$ in Eq. (12). □

APPENDIX III

Because the $\ell^2$ norm of the error is proportional to the $\ell^2$ norm of the solution, as indicated by Theorem 1, we prove Theorem 2 by showing that the $\ell^2$ norm of the solution is decaying according to Eq. (15).

The solution to the first iteration of the outer loop ($l$=1) is $\mathbf{x}^{(1)} = \mathbf{x}^*$. The first outer loop produces an estimate $\delta\mathbf{x}$, with an error $|\mathbf{x}^* - \delta\mathbf{x}^{(1)}|_2$ bounded as in Eq. (12). The second iteration of the outer loop solves the linear system $\mathbf{r}^{(1)} = \mathbf{A}\delta\mathbf{x}$, where $\mathbf{r}^{(1)} = \mathbf{y} - \mathbf{A}\delta\mathbf{x}^{(1)}$, and its solution is $\mathbf{x}^{(2)} = \mathbf{x}^* - \delta\mathbf{x}^{(2)}$, whose $\ell^2$ norm is bounded by $|\mathbf{x}^{(2)}|_2 \leq \theta|\mathbf{x}^*|_2$. Using Eq. (12) in Theorem 1, the $\ell^2$ norm of the error between the solution from the second set of Richardson iterations $\delta\mathbf{x}^{(2)}$ and its solution $\mathbf{x}^{(2)}$ is bounded by,

$$|\mathbf{x}^{(2)} - \delta\mathbf{x}^{(2)}|_2 \leq \theta|\mathbf{x}^{(2)}| \leq \theta^2|\mathbf{x}^*|_2. \quad (32)$$

Rewriting $|\mathbf{x}^{(2)} - \delta\mathbf{x}^{(2)}|_2$ as $|\mathbf{x}^* - \delta\mathbf{x}^{(1)} - \delta\mathbf{x}^{(2)}|_2 = |\mathbf{x}^* - \mathbf{x}^{(1)}|_2$, the accumulated solution $\mathbf{x}^{(2)} = \delta\mathbf{x}^{(1)} + \delta\mathbf{x}^{(2)}$ has an error $|\mathbf{x}^* - \mathbf{x}^{(1)}|_2 \leq \theta^2|\mathbf{x}^*|_2$, which is at least $\theta$ times smaller than the solution from the first iteration of the outer loop $\delta\mathbf{x}^{(1)}$. After $M$ iterations of outer loop, $|\mathbf{x}^* - \mathbf{x}^{(M-1)}|_2 \leq \theta^M|\mathbf{x}^*|_2$. □

ACKNOWLEDGMENT

The authors would like to thank Dr. Stephen Becker (Department of Applied Mathematics, University of Colorado Boulder) for helpful discussions.